\documentclass[a4paper,11pt,headsepline,twoside,leqno]{scrartcl}
\usepackage[T1]{fontenc}
\usepackage[left=25mm, right=25mm, top=25mm, bottom=27mm]{geometry}
\usepackage{amssymb, amsmath, amsthm}
\usepackage{tikz}
\usepackage{thmtools}
\usepackage{authblk}
\usepackage{helvet}
\usepackage[square,authoryear]{natbib}
\usepackage[labelsep=period,labelfont=bf,font=small]{caption}
\usepackage{subcaption}
\usepackage{graphicx}
\usepackage{mathtools}
\usepackage[hang]{footmisc}
\usepackage{float}
\usepackage{fancyhdr}
\usepackage{setspace}
\usepackage{afterpage}
\usepackage[breaklinks]{hyperref}

\setcapindent{0pt}

\setkomafont{sectioning}{\normalfont\normalcolor\bfseries}

\floatstyle{plaintop}
\restylefloat{table}

\title{\LARGE{Translation matrix elements for spherical Gauss-Laguerre basis functions}\vspace*{11pt}}

\author[1]{\large{J\"urgen Prestin}\vspace*{11pt}}
\author[,2]{\large{Christian W\"ulker}\thanks{\!Email:\ christian.wuelker@jhu.edu}}

\affil[1]{\normalsize{\textit{Institute of Mathematics, University of L\"ubeck, Germany}}\vspace*{11pt}}

\affil[2]{\normalsize{\textit{Department of Mechanical Engineering,} 
\vspace*{2pt}\authorcr\textit{Johns Hopkins University, Baltimore, MD, USA}}\vspace*{11pt}}

\date{\normalsize{May 20, 2018}}

\newtheoremstyle{definition}%
  {5pt}
  {5pt}
  {}
  {}
  {\bfseries}
  {\textbf{.}}
  {.5em}
  {\thmname{#1} \thmnumber{#2}}
  
\theoremstyle{definition} 
\newtheorem{definition}{Definition}
\newtheorem{problem}[definition]{Problem}

\newtheoremstyle{definitioncited}%
  {5pt}
  {5pt}
  {}
  {}
  {\bfseries}
  {\textbf{.}}
  {.5em}
  {\thmname{#1} \thmnumber{#2} \thmnote{\normalfont#3}}
  
\theoremstyle{definitioncited}  
\newtheorem{definitioncited}[definition]{Definition} 

\newtheoremstyle{theoremcited}%
  {5pt}
  {5pt}
  {\itshape}
  {}
  {\bfseries}
  {\textbf{.}}
  {.5em}
  {\thmname{#1} \thmnumber{#2} \thmnote{\normalfont#3}}

\theoremstyle{theoremcited}  
\newtheorem{maintheoremcited}[definition]{Main theorem}
\newtheorem{theoremcited}[definition]{Theorem}

\newtheorem{lemmacited}[definition]{Lemma}

\newtheoremstyle{theorem}%
  {5pt}
  {5pt}
  {\itshape}
  {}
  {\bfseries}
  {\textbf{.}}
  {.5em}
  {\thmname{#1} \thmnumber{#2}}

\theoremstyle{theorem}  

\newtheorem{theorem}[definition]{Theorem}
\newtheorem{lemma}[definition]{Lemma}

\newtheoremstyle{remark}%
  {5pt}
  {5pt}
  {\normalfont}
  {}
  {\itshape}
  {.}
  {.5em}
  {\thmname{#1} \thmnumber{#2}}

\theoremstyle{remark}
\newtheorem{remark}[definition]{Remark}

\numberwithin{definition}{section}
\numberwithin{definitioncited}{section}
\numberwithin{corollary}{section}
\numberwithin{corollarycited}{section}
\numberwithin{maintheorem}{section}
\numberwithin{maintheoremcited}{section}
\numberwithin{theorem}{section}
\numberwithin{theoremcited}{section}
\numberwithin{lemma}{section}
\numberwithin{lemmacited}{section}
\numberwithin{problem}{section}
\numberwithin{remark}{section}

\numberwithin{equation}{section}

\declaretheoremstyle[%
  spaceabove=-1pt,%
  spacebelow=5pt,%
  headfont=\normalfont\itshape,%
  postheadspace=5pt,%
  qed=\qedsymbol%
]{mystyle} 
\declaretheorem[name={Proof},style=mystyle,unnumbered,
]{prf}
\usetikzlibrary{calc,fadings,decorations.pathreplacing}

\newcommand\pgfmathsinandcos[3]{%
  \pgfmathsetmacro#1{sin(#3)}%
  \pgfmathsetmacro#2{cos(#3)}%
}
\newcommand\LongitudePlane[3][current plane]{%
\pgfmathsinandcos\sinEl\cosEl{#2} 
\pgfmathsinandcos\sint\cost{#3} 
\tikzset{#1/.estyle={cm={\cost,\sint*\sinEl,0,\cosEl,(0,0)}}}
}
\newcommand\LatitudePlane[3][current plane]{%
  \pgfmathsinandcos\sinEl\cosEl{#2} 
  \pgfmathsinandcos\sint\cost{#3} 
  \pgfmathsetmacro\yshift{\cosEl*\sint}
  \tikzset{#1/.style={cm={\cost,0,0,\cost*\sinEl,(0,\yshift)}}} %
}
\newcommand\DrawLongitudeCircle[2][1]{
\LongitudePlane{\angEl}{#2}
\tikzset{current plane/.estyle={cm={\cost,\sint*\sinEl,0,\cosEl,(0,0)},scale=#1}}
\pgfmathsetmacro\angVis{atan(sin(#2)*cos(\angEl)/sin(\angEl))} %
\draw[current plane,thin,black] (\angVis:1) arc (\angVis:\angVis+180:1);
\draw[current plane,thin,dashed] (\angVis-180:1) arc (\angVis-180:\angVis:1);
}
\newcommand\DrawLatitudeCircle[2][1]{
  \LatitudePlane{\angEl}{#2}
  \tikzset{current plane/.prefix style={scale=#1}}
  \pgfmathsetmacro\sinVis{sin(#2)/cos(#2)*sin(\angEl)/cos(\angEl)}
  \pgfmathsetmacro\angVis{asin(min(1,max(\sinVis,-1)))}
  \draw[current plane] (\angVis:1) arc (\angVis:-\angVis-180:1);
  \draw[current plane,dashed] (180-\angVis:1) arc (180-\angVis:\angVis:1);
}

\tikzset{%
  >=latex, 
  inner sep=0pt,%
  outer sep=2pt,%
  mark coordinate/.style={inner sep=0pt,outer sep=0pt,minimum size=3pt,
    fill=black,circle}%
}

\setkomafont{sectioning}{\rmfamily\bfseries}
\setkomafont{section}{\large\rmfamily\bfseries}
\setkomafont{subsection}{\rmfamily\bfseries}
\setkomafont{subsubsection}{\rmfamily\bfseries}

\begin{document}

\maketitle
\thispagestyle{empty}

\begin{abstract}
\begin{center} 
\textbf{Abstract}\vspace*{11pt}
\end{center}
\renewcommand{\baselinestretch}{1.4}\normalsize
Spherical Gauss-Laguerre (SGL) basis functions, \textit{i.e.}, normalized functions of the type $L_{n-l-1}^{(l + 1/2)}(r^2) r^{l} Y_{lm}(\vartheta,\varphi)$, $|m| \leq l < n \in \mathbb{N}$, 
constitute an orthonormal polynomial basis of the space $L^{2}$ on $\mathbb{R}^{3}$ with radial Gaussian weight $\exp(-r^{2})$. We have recently described reliable fast Fourier transforms for the SGL basis functions. The main application of the SGL basis functions and our fast algorithms is in solving certain three-dimensional rigid matching problems, where the center is prioritized over the periphery. For this purpose, so-called SGL translation matrix elements are required, which describe the spectral behavior of the SGL basis functions under translations. In this paper, we derive a closed-form expression of these translation matrix elements, allowing for a direct computation of these quantities 
in practice. 
\\[11pt]
\textit{Key words and phrases:}\vspace*{-4pt} Spherical Gauss-Laguerre (SGL) basis functions, translation, three-dimensional rigid matching, computational harmonic analysis
\end{abstract}

\setstretch{1.125}

\newpage
\section{Introduction}
\emph{Spherical Gauss-Laguerre (SGL) basis functions} (Def.\:\ref{def:sgl}) constitute an orthonormal polynomial basis of the space $L^{2}$ on $\mathbb{R}^{3}$ equipped with the radial Gaussian weight function $\exp(-r^{2})$. We have recently described reliable fast Fourier transforms for the SGL basis functions \citep{prestin_wuelker}. These algorithms allow for a fast computation of SGL Fourier coefficients for spectral analysis. A main application of the SGL basis functions and our SGL Fourier transforms is the fast solution of certain three-dimensional rigid matching problems, where the center is prioritized over the periphery (Prob.\:\ref{prob:1}). For this purpose, so-called \textit{SGL translation matrix elements} are required. These elements describe the spectral behavior of the SGL basis functions under translations in $\mathbb{R}^{3}$. In this paper, we derive a closed-form expression of these translation matrix elements (Thm.\:\ref{thm:sgl_translation_matrix_elements}). This allows for a direct computation of the SGL translation matrix elements when required in practice. 

In the derivation of the closed-form expression of the SGL translation matrix elements, we make use of the fact that the SGL basis functions are akin to so-called \textit{Gaussian-type orbital (GTO) basis functions} in the unweighted space $L^{2}(\mathbb{R}^{3})$, which themselves carry an exponential radial decay factor $\exp(-r^{2}/2)$ \citep{ritchie_kemp}. These functions are nowadays extensively made use of in biomolecular recognition simulation such as protein-protein docking (see, \textit{e.g.}, \citep{ritchie_kemp} or \citep[Sec.\:2.1]{ritchie_kozakov_vajda}). \citet{ritchie} has already established a closed-form expression of the GTO translation matrix elements. In principle, we can apply the same approach to the SGL case. However, as a key tool, \citeauthor{ritchie} introduces a special class of Hankel transforms, referred to as the \emph{spherical Bessel transform} \citep[Sec.\:2.3]{ritchie}. The GTO basis functions are eigenfunctions of these transforms, a fact that crucially simplifies the derivation of the GTO translation matrix elements. This is not the case with the SGL basis functions. In fact, the spherical Bessel transform of \citeauthor{ritchie} is not even well-defined in the SGL case, as the underlying integrals diverge. For this reason, we introduce a \textit{weighted spherical Bessel transform} (Def.\:\ref{def:weighted_spherical_Bessel_transform}) by adding a Gauss-Weierstrass convergence factor to the unweighted spherical Bessel transform used by \citeauthor{ritchie}. While this ensures convergence of the corresponding integrals, it also causes additional technical difficulties we must solve. 

The remainder of this paper is organized as follows: In Section \ref{sec:background}, we give a brief overview on SGL basis functions. We illustrate their application, and show how our fast SGL Fourier transforms can be used in this context. In Section \ref{sec:theory}, we present the theoretical tools to derive the closed-form expression for the SGL translation matrix elements. The derivation itself is presented in Section \ref{sec:proof}.

\section{Background}\label{sec:background}
Let $\|\cdot\|_{2}$ denote the standard Euclidean norm on $\mathbb{R}^{3}$. We consider the weighted $L^{2}$ space
\begin{equation*}
H \hspace*{2pt}\coloneqq\hspace*{2pt} \left\lbrace f : \mathbb{R}^{3\!} \to \mathbb{C} 
\hspace*{1pt}:\hspace*{1pt} f~\textnormal{Lebesgue measurable and}\hspace*{1pt}
\int_{\mathbb{R}^{3}} |f(\boldsymbol{x})|^{2} \hspace*{1pt} \mathrm{e}^{-\|\boldsymbol{x}\|_{2}^{2}} \hspace*{1pt} \mathrm{d}\boldsymbol{x} < \infty 
\right\rbrace\!,
\end{equation*}
endowed with the inner product
\begin{equation*}\label{eq:inner_product}
\langle f, g \rangle_{H} \hspace*{2pt}\coloneqq\hspace*{2pt} \int_{\mathbb{R}^{3}} f(\boldsymbol{x}) \hspace*{1pt} \overline{g(\boldsymbol{x})} \hspace*{1pt} \mathrm{e}^{-\|\boldsymbol{x}\|_{2}^{2}} \hspace*{1pt} \mathrm{d}\boldsymbol{x}, ~~~~ f, g \in H,
\end{equation*}
and induced norm $\|\cdot\|_{H} \coloneqq \sqrt{\langle\cdot,\cdot\rangle_{H}}$.

SGL basis functions, evolutionarily related to the GTO basis functions in the unweighted space $L^{2}(\mathbb{R}^{3})$ introduced by \citet{ritchie_kemp}, are orthogonal polynomials in the weighted $L^{2}$ space $H$. They arise from a particular construction approach in spherical coordinates. These are defined as \emph{radius} $r \in [0,\infty)$, \emph{polar angle} $\vartheta \in [0,\pi]$, and \emph{azimuthal angle} $\varphi \in [0,2\pi)$, being connected to Cartesian coordinates $x$, $y$, and $z$ via (cf.\ Fig.\:\ref{fig:kugelkoordinaten})
\begin{align*}
x \hspace*{2pt}&=\hspace*{2pt} r \hspace*{1pt} \sin\vartheta \hspace*{1pt} \cos\varphi,\\
y \hspace*{2pt}&=\hspace*{2pt} r \hspace*{1pt} \sin\vartheta \hspace*{1pt} \sin\varphi,\\
z \hspace*{2pt}&=\hspace*{2pt} r \hspace*{1pt} \cos\vartheta.
\end{align*}
In this paper, we write $\boldsymbol{x} = [r,\vartheta,\varphi]$ when $r$, $\vartheta$, and $\varphi$ are spherical coordinates of $\boldsymbol{x}$ (potential ambiguity is not a problem). By \hspace*{1pt}$\mathbb{S}^{2} \coloneqq \lbrace \boldsymbol{x} \in \mathbb{R}^{3} : \|\boldsymbol{x}\|_{2} = 1 \rbrace$ we denote the two-dimensional unit sphere. The SGL basis functions are now defined as follows.

\begin{definitioncited}[\textbf{(SGL basis functions)}]\label{def:sgl}
The SGL basis function of orders $n \in \mathbb{N}$, $l \in \{0,\dots,$ $n-1\}$, and $m \in \{-l,\dots,l\}$ is defined as
\begin{equation*}\label{eq:sgl}
H_{nlm} : \mathbb{R}^{3} \to \mathbb{C}, ~~~~ H_{nlm} (r,\vartheta,\varphi) \hspace*{0pt}\coloneqq\hspace*{0pt} N_{nl} \hspace*{1pt} R_{nl} (r) \hspace*{1pt} Y_{lm}(\vartheta,\varphi),
\end{equation*}
where 
\begin{equation*}
N_{nl} \hspace*{2pt}\coloneqq\hspace*{2pt} \sqrt{\frac{2(n-l-1)!}{\Gamma(n + 1/2)}} 
\end{equation*}
is a normalization constant, $Y_{lm}$ is the spherical harmonic of degree $l$ and order $m$ \citep[Sec.\:1.6.2]{dai_xu}, while the radial part $R_{nl}$ is defined 
as
\begin{equation*}
R_{nl}(r) \hspace*{2pt}\coloneqq\hspace*{2pt} L_{n-l-1}^{(l+1/2)\!}(r^{2}) \hspace*{1pt} r^{l},
\end{equation*}
$L_{n - l - 1}^{(l + 1/2)}$ being a generalized Laguerre polynomial \citep[Sec.\:6.2]{andrews_askey_roy}.
\end{definitioncited}

\begin{theoremcited}[\textbf{\citep[Cor.\:1.3]{prestin_wuelker}}]
The SGL basis functions constitute an orthonormal polynomial basis of the space \hspace*{1pt}$H$. Thus, every function \hspace*{1pt}$f \in H$ can be approximated arbitraly well with respect to \hspace*{1pt}$\|\cdot\|_{H}$ by linear combinations of the SGL basis functions.
\end{theoremcited}

For a translation vector $\boldsymbol{t} \in \mathbb{R}^{3}$, we define the \emph{translation operator} $T(\boldsymbol{t})$ via
\begin{equation*}
T(\boldsymbol{t})[f](\boldsymbol{x}) \hspace*{2pt}\coloneqq\hspace*{2pt} f(\boldsymbol{x}-\boldsymbol{t}).
\end{equation*}
Further, for a given rotation $R \in \textnormal{SO}(3)$ in $\mathbb{R}^{3}$, we define the \emph{rotation operator} $R$ via
\begin{equation*}
Rf(\boldsymbol{x}) \hspace*{2pt}\coloneqq\hspace*{2pt} f(R^{-1}\boldsymbol{x}).
\end{equation*}
The main application of the SGL basis functions is the fast solution of the following three-dimensional rigid matching problem, where due to the present weight function, the center is prioritized over the periphery:

\begin{figure}[t]
\begin{center}
\begin{tikzpicture}
\input{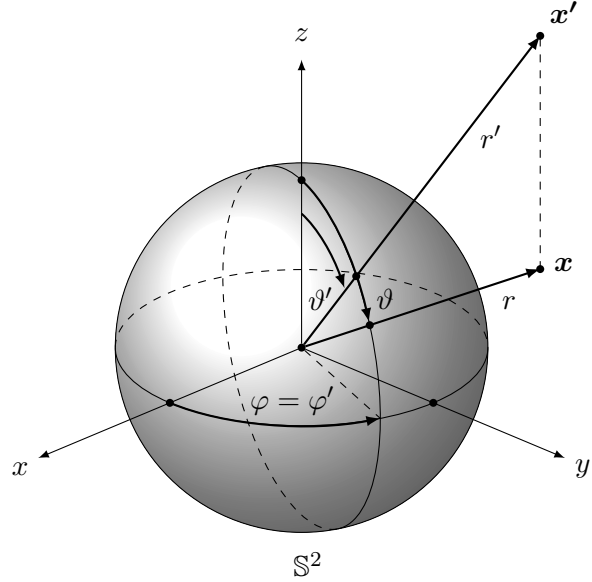}
\end{tikzpicture}
\end{center}
\caption{Spherical coordinates in $\mathbb{R}^{3}$. The point $\boldsymbol{x'}$ with spherical coordinates $r'$, $\vartheta'$, and $\varphi'$ arises from the point $\boldsymbol{x}$ with spherical coordinates $r$, $\vartheta$, and $\varphi$ by a translation along the $z$ axis.}\label{fig:kugelkoordinaten}
\end{figure}

\begin{problem}\label{prob:1}
For admissible functions $f$ and $g$, maximize 
with respect to $R \in \textnormal{SO}(3)$ and $\boldsymbol{t} \in \mathbb{R}^{3}$
the absolute value of the weighted overlap integral
\begin{equation}\label{eq:overlap_integral_1}
I(R,\boldsymbol{t}) \hspace*{2pt}\coloneqq\hspace*{2pt} \int_{\mathbb{R}^{3}}\! T(\boldsymbol{t}) [R \hspace*{1pt} f](\boldsymbol{x}) \hspace*{1pt} \overline{g(\boldsymbol{x})} \hspace*{1pt} \mathrm{e}^{-\|\boldsymbol{x}\|_{2}^{2}} \hspace*{1pt} \mathrm{d}\boldsymbol{x}.
\end{equation}
%
\end{problem}


Problem \ref{prob:1} can be tackled in the following way: Assume that $f, g \in H$ are \emph{bandlimited functions} with \emph{bandwidth} $B \in \mathbb{N}$, \textit{i.e.}, $\langle f, H_{nlm} \rangle_{H} = \langle g, H_{nlm} \rangle_{H} = 0$ for $n > B$. This means that
\begin{align*}
f \hspace*{2pt}&=\hspace*{2pt} \sum\nolimits_{|m| \leq l < n \leq B} \hat{f}_{nlm} \hspace*{1pt} H_{nlm},\\[4pt]
g \hspace*{2pt}&=\hspace*{2pt} \sum\nolimits_{|m| \leq l < n \leq B} \hat{g}_{nlm} \hspace*{1pt} H_{nlm},
\end{align*}
where the sums run over the index ranges for the SGL basis functions of order at most $n \leq B$, and $\hat{f}_{nlm} \coloneqq \langle f, H_{nlm} \rangle_{H}$ and $\hat{g}_{nlm} \coloneqq \langle g, H_{nlm} \rangle_{H}$ are the \emph{SGL Fourier coefficients} of $f$ and $g$, respectively. The overlap integral in \eqref{eq:overlap_integral_1} thus attains the form
\begin{equation}\label{eq:overlap_integral_2}
I(R,\boldsymbol{t})
\hspace*{2pt}=\hspace*{2pt} \sum_{|m| \leq l < n \leq B} \hat{f}_{nlm} \sum_{|m'| \leq l'\! < n' \leq B} \overline{\hat{g}_{n'l'm'}} \underbrace{\int_{\mathbb{R}^{3}}\! T(\boldsymbol{t}) [R \hspace*{1pt} H_{nlm}](\boldsymbol{x}) \hspace*{1pt} \overline{H_{n'l'm'}(\boldsymbol{x})} \hspace*{1pt} \mathrm{e}^{-\|\boldsymbol{x}\|_{2}^{2}} \hspace*{1pt} \mathrm{d}\boldsymbol{x}}_{\langle T(\boldsymbol{t}) R \hspace*{1pt} H_{nlm}, \hspace*{1pt} H_{n'l'm'} \rangle_{H}}\!.
\end{equation}
This approach has the major advantage that the SGL Fourier coefficients of $f$ and $g$ need to be computed only once for evaluating the overlap integral $I$ for several different rotations and translations. For exactly this purpose, we have recently described fast Fourier transforms for the SGL basis functions\:--\:see \citep{prestin_wuelker} for gridded data, and \citep[Chap.\:3]{wuelker} for scattered data. The fast algorithms for gridded data are based on an exact quadrature formula with $(2B)^{3}$ sampling points in $\mathbb{R}^{3}$  for the $B(B+1)(2B+1)/6$ potentially nonzero SGL Fourier coefficients of bandlimited functions with bandlimit $B$. The asymptotic complexity of these algorithms is $\mathcal{O}(B^{4})$ or even $\mathcal{O}(B^{3} \log^{2} B)$, instead of the naive complexity of $\mathcal{O}(B^{6})$\footnote{Our C++ implementation of these fast algorithms is available from \url{https://github.com/cwuelker/SGLPack}}. Furthermore, the integral on the right-hand side of \eqref{eq:overlap_integral_2} 
is now independent of the particular functions $f$ and $g$ considered, and can hence be computed independently. This integral combines the spectral behavior of the SGL basis functions under rotations and translations.

The spectral behavior of the SGL basis functions under rotations is relatively simple to describe; it is directly inherited from the spherical harmonics:
\begin{theorem}\label{thm:sgl_rotation}
For \hspace*{1pt}$R \in \textnormal{SO}(3)$, we have that
\begin{equation*}
\langle R \hspace*{1pt} H_{nlm}, H_{n'l'm'} \rangle_{H} \hspace*{2pt}=\hspace*{2pt} \delta_{nn'} \hspace*{1pt} \delta_{ll'} \hspace*{1pt} D_{mm'}^{(l)}(R),
\end{equation*}
where \hspace*{1pt}$D_{mm'}^{(l)}$ are the well-known \emph{Wigner-$D$ functions} \citep[Sec.\:10.5\:ff.]{nikiforov_uvarov}.
\end{theorem}
%

The 
behavior 
under translations is somewhat more complicated. With regard to Problem \ref{prob:1}, it is sufficiently described by the \emph{SGL translation matrix elements}
\begin{equation*}
T_{nn'll'}^{(mm')}(\nu) \hspace*{2pt}\coloneqq\hspace*{2pt} \langle T(\nu \textbf{e}_{z}) H_{nlm}, H_{n'l'm'} \rangle_{H}, ~~~~ \nu > 0,
\end{equation*}
where $\textbf{e}_{z}$ denotes the unit vector along the positive $z$ axis. As explained below, due to the rotation invariance of the Lebesgue measure on $\mathbb{R}^{3}$, it is sufficient to consider solely translations along the $z$ axis in this context. This also has the advantage that the translation operator $T(\nu \textbf{e}_{z})$ has no impact on the azimuthal angle $\varphi$ (cf.\ Fig.\:\ref{fig:kugelkoordinaten}), a fact that simplifies the derivation of the closed-form expression for the SGL translation matrix elements and the closed-form expression itself significantly.



Let us assume now that we want to evaluate the Integral $I$ in \eqref {eq:overlap_integral_1}  for $M$ pairs $(R_{i},\boldsymbol{t}_{i}) \in \textnormal{SO}(3) \times \mathbb{R}^{3}$, $i = 0, \dots, M-1$, in order to empirically determine its maximum absolute value, that is, to solve Problem \ref{prob:1}. We propose the following strategy: In a first step, we compute the SGL Fourier coefficients $\hat{f}_{nlm}$ and $\hat{g}_{nlm}$ of $f$ and $g$, respectively, using our fast SGL Fourier transforms. Then, for every $i = 0, \dots, M - 1$, we choose $\tilde{R}_{i} \in \textnormal{SO}(3) $ such that $\tilde{R}_{i} \hspace*{1pt} \boldsymbol{t}_{i}$ points in the positive direction of the $z$ axis. With the aid of the rotation invariance of both the Lebesgue measure and the weight function $\exp(-r^{2})$ on $\mathbb{R}^{3}$, and by making use of 
Theorem \ref {thm:sgl_rotation}, 
we get
\begin{align}\label{eq:inner_product_split}
\langle T(\boldsymbol{t}_{i}) \hspace*{1pt} R_{i} \hspace*{1pt} H_{nlm}, H_{n'l'm'} \rangle_{H}
\hspace*{2pt}&=\hspace*{2pt} \langle \tilde{R}_{i} \hspace*{1pt} T(\boldsymbol{t}_{i}) \hspace*{1pt} R_{i} \hspace*{1pt} H_{nlm}, \tilde{R}_{i} \hspace*{1pt} H_{n'l'm'} \rangle_{H}\\[6pt]\notag
\hspace*{2pt}&=\hspace*{2pt} \langle T(\tilde{R}_{i} \hspace*{1pt} \boldsymbol{t}_{i}) \hspace*{1pt} \tilde{R}_{i} \hspace*{1pt} R_{i} \hspace*{1pt} H_{nlm}, \tilde{R}_{i} \hspace*{1pt} H_{n'l'm'} \rangle_{H}\\[2pt]\notag
\hspace*{2pt}&=\hspace*{2pt} \sum_{|m''| \leq l} D^{(l)}_{mm''}(\tilde{R}_{i} \hspace*{1pt} R_{i}) \sum_{|m'''| \leq l'} \overline{D^{(l')}_{m'm'''}(\tilde{R}_{i})} \hspace*{2pt} T_{nn'll'}^{(m''m''')}(\|\boldsymbol{t}_{i}\|_{2}).
\end{align}
We will later see that generally $T_{nn'll'}^{(mm')} = \delta_{mm'} \hspace*{1pt} T_{nn'll'}^{(|m|)}$, \textit{i.e.}, the SGL translation matrix elements vanish unless $m = m'$, and they do not depend on the sign of $m$. Now the values $I(R_{i},\boldsymbol{t}_{i})$, $i = 0, \dots, M - 1$, can be computed using \eqref{eq:overlap_integral_2} and \eqref{eq:inner_product_split}, together with the closed-form expression for the SGL translation matrix elements in Theorem \ref{thm:sgl_translation_matrix_elements}.

%
%
%

\section{Theory}\label{sec:theory}
The main idea in the derivation of a closed-form expression for the GTO translation matrix elements of \citet{ritchie} is to establish a connection between the GTO basis functions and certain products of spherical Bessel functions and spherical harmonics. The behavior of the arising functions under translations along the $z$ axis is known, and this knowledge can be transferred to the GTO (and thus the SGL) case. The spherical Bessel functions are defined as follows:

\begin{definitioncited}[\textbf{(Spherical Bessel functions)}]
For $n \in \mathbb{N}_{0}$, the \emph{spherical Bessel function (of first kind)} is defined as \citep[Sec.\:10.1.1]{abramowitz_stegun}
\begin{equation*}
j_{n} : (0,\infty) \to \mathbb{R}, ~~~~ j_{n}(\xi) \hspace*{2pt}\coloneqq\hspace*{2pt} \sqrt{\frac{\pi}{2 \xi}} \hspace*{1pt} J_{n+1/2}(\xi),
\end{equation*}
where $J_{n+1/2}$ denotes the Bessel function (of first kind) of fractional order $n+1/2$ (see \citep{watson} for a detailed introduction).
\end{definitioncited}

\begin{remark}
Via the limit \citep[Eq.\:10.1.4]{abramowitz_stegun}
\begin{equation*}
\lim_{\xi \to 0+} \hspace*{0pt} j_{n}(\xi) \hspace*{2pt}=\hspace*{2pt}
\begin{cases}
1, \!\!\!\!\!\!&\textnormal{ if } n = 0,\\
0  \!\!\!\!\!\!&\textnormal{ otherwise},
\end{cases}
\end{equation*}
the spherical Bessel functions can be extended continuously to $[0,\infty)$. 
\end{remark}

The behavior of the product functions $j_{l}(r) Y_{lm}(\vartheta,\varphi)$ under translations along the $z$ axis is described by the important next lemma. Recall that translations along the $z$ axis do not affect the azimuthal angle $\varphi$, as mentioned above.

\begin{lemmacited}[\textbf{(Spherical Bessel addition theorem)}]\label{lemma:Bessel_addition_theorem}
Let \hspace*{1pt}$|m| \leq l \in \mathbb{N}_{0}$ and \hspace*{1pt}$\boldsymbol{x} = [r,\vartheta,\varphi]$. Let further \hspace*{1pt}$\nu > 0$ be given and set \hspace*{1pt}$\boldsymbol{x}'\! \coloneqq \boldsymbol{x} - \nu \mathbf{e}_{z} = [r'\!,\vartheta'\!,\varphi]$. Then
\begin{equation}\label{eq:Bessel_addition_theorem}
j_{l}(\beta r) \hspace*{1pt} Y_{lm}(\vartheta,\varphi) \hspace*{2pt}=\hspace*{2pt} \sum_{l'=0}^{\infty} \hspace*{1pt} \sum_{k = |l-l'|}^{l + l'\!} A_{k}^{(ll'm)} \hspace*{1pt} j_{k}(\beta \nu) \hspace*{1pt} j_{l'}(\beta r') \hspace*{1pt} Y_{l'm}(\vartheta'\!,\varphi), ~~~~ \beta \in (0,\infty),
\end{equation}
where 
\begin{equation}\label{eq:A_coefficients}
A_{k}^{(ll'm)} \hspace*{1pt}\coloneqq\hspace*{2pt} (-1)^{(k-l+l')/2 + m} \hspace*{1pt} \sqrt{(2l+1)(2l'+1)} \hspace*{1pt} (2k+1)
\begin{pmatrix}
l \!&\! l' \!&\! k \\
0 \!&\! 0  \!&\! 0
\end{pmatrix}\!\!
\begin{pmatrix}
l \!&\! l' \!&\! k \\
m \!&\! -m \!&\! 0
\end{pmatrix}\!.
\end{equation}
The series in \eqref{eq:Bessel_addition_theorem} converges uniformly with respect to \hspace*{1pt}$[\vartheta'\!,\varphi] \in \mathbb{S}^{2}$.
\end{lemmacited}

\begin{remark}
On the right-hand side of \eqref{eq:A_coefficients}, we find the well-known \emph{Wigner-$3j$ symbols} (see \citep[Sec.\:3.12]{biedenharn_louck} for a detailed introduction). Furthermore, it is easy to verify that the coefficients $A_{k}^{(ll'm)}$ 
are independent of the sign of $m$, that is $A_{k}^{(ll'm)} = A_{k}^{(ll'|m|)}$.
\end{remark}

\noindent\textit{Proof of Lemma \ref{lemma:Bessel_addition_theorem}.}
We follow the outline of \citet[Appx.\:A]{ritchie}, while clarifying some technical details not considered there. Combining Bauer's Bessel addition theorem \citep[Sec.\:11.5, Eq.\:1]{watson} and the addition theorem of the spherical harmonics \citep[Eq.\:1.6.7]{dai_xu} 
%
gives the \emph{plane wave expansion}
\begin{equation}\label{eq:plane_wave_expansion}
\mathrm{e}^{\mathrm{i} \langle \boldsymbol{x},\boldsymbol{k} \rangle_{2}} \hspace*{2pt}=\hspace*{2pt} 4\pi \hspace*{1pt} \sum\nolimits_{|m| \leq l \in \mathbb{N}_{0}} \mathrm{i}^{l} j_{l}(\beta r) \hspace*{1pt} Y_{lm}(\vartheta,\varphi) \hspace*{1pt} \overline{Y_{lm}(\theta,\phi)}, ~~~~
%
\boldsymbol{k} = [\beta,\theta,\phi],
\end{equation}
$\boldsymbol{k}$ being referred to as the \emph{wave vector}. 
The series on the right-hand side of \eqref{eq:plane_wave_expansion} is absolutely convergent, which can be shown by estimating (cf.\ \citep[P.\:53]{watson})
\begin{equation}\label{eq:Bessel_estimate}
|J_{n+1/2}(\xi)| \hspace*{2pt}\leq\hspace*{2pt}  \frac{2 \hspace*{1pt} (\xi/2)^{n+1/2}}{\sqrt{\pi} \hspace*{1pt} n!}, ~~~~ n \in \mathbb{N}_{0}, ~ \xi \in [0,\infty),
\end{equation}
and \citep[Eq.\:3.1.4]{freeden_gervens_schreiner}
\begin{equation}\label{eq:Y_estimate}
|Y_{lm}| \hspace*{2pt}\leq\hspace*{2pt} \sqrt{\frac{2l+1}{4\pi}}, ~~~~ |m| \leq l \in \mathbb{N}_{0}.
\end{equation}

Let now $\nu \mathbf{e}_{z} = [\nu,0,0]$. By \eqref{eq:plane_wave_expansion} and using 
%
$\mathrm{e}^{\mathrm{i} \langle \boldsymbol{x},\boldsymbol{k} \rangle_{2}} 
= \mathrm{e}^{\mathrm{i} \langle \boldsymbol{x}'\!, \boldsymbol{k} \rangle_{2}} \hspace*{1pt} \mathrm{e}^{\mathrm{i} \langle \nu \mathbf{e}_{z},\boldsymbol{k} \rangle_{2}}$,
%
we get 
\begin{align*}
\sum\nolimits_{|m| \leq l \in \mathbb{N}_{0}} \mathrm{i}^{l} j_{l}(\beta r) \hspace*{1pt} Y_{lm}(\vartheta,\varphi) \hspace*{1pt} \overline{Y_{lm}(\theta,\phi)}
\hspace*{2pt}&=\hspace*{2pt} 4\pi \sum\nolimits_{|m'| \leq l' \in \mathbb{N}_{0}} \sum\nolimits_{|m''| \leq l'' \in \mathbb{N}_{0}} \hspace*{1pt}\mathrm{i}^{l'\!+l''\!} j_{l'}(\beta r') j_{l''}(\beta \nu)\\\notag
&\hspace*{14pt} \times Y_{l'm'}(\vartheta'\!,\varphi) \hspace*{1pt} Y_{l''m''}(0,0) \hspace*{1pt} \overline{Y_{l'm'}(\theta,\phi) \hspace*{1pt} Y_{l''m''}(\theta,\phi)}.
\end{align*}
%
Both sides of the equation converge absolutely and uniformly with respect to $[\theta,\phi] \in \mathbb{S}^{2}$, which follows again from \eqref{eq:Bessel_estimate} and \eqref{eq:Y_estimate}, and applying Weierstrass' test. 
Thus, using the orthonormality and the identity $\overline{Y_{lm}} = (-1)^{m}Y_{l,-m}$ of the spherical harmonics, we find that 
\begin{align}\label{eq:spherical_bessel_addition_theorem_zwischenschritt_2}
j_{l}(\beta r) \hspace*{1pt} Y_{lm}(\vartheta,\varphi)
\hspace*{2pt}&=\hspace*{2pt} 4\pi \hspace*{0.5pt}\sum\nolimits_{|m'| \leq l' \in \mathbb{N}_{0}} \sum\nolimits_{|m''| \leq l'' \in \mathbb{N}_{0}} \mathrm{i}^{-l+l'\!+l''\!} j_{l'}(\beta r') \hspace*{1pt} j_{l''}(\beta \nu) \hspace*{1pt} Y_{l'm'}(\vartheta'\!,\varphi) \hspace*{1pt} Y_{l''m''}(0,0)\\\notag
&\hspace*{14pt} \times (-1)^{m'\!+m''}\!\! \int_{0}^{\pi} \!\!\int_{0}^{2\pi}\!\! Y_{lm}(\theta,\phi) \hspace*{1pt} Y_{l'\!,-m'}(\theta,\phi) \hspace*{1pt} Y_{l''\!,-m''}(\theta,\phi) \hspace*{2pt} \mathrm{d}\phi \hspace*{1pt} \sin\theta \hspace*{2pt} \mathrm{d}\theta.
\end{align}
Similarly, due to the orthonormality of the trigonometric monomials $\exp(\mathrm{i}m\varphi)$, we can discard the summation over $m'$ on the right-hand side, leaving only the summands with $m'\!=m$.
For the double integral on the right-hand side, we then use Gaunt's formula for integrals over three spherical harmonics \citep[Eq.\:3.192]{biedenharn_louck} to obtain
\begin{align}\label{eq:Gaunt}
\int_{0}^{\pi} \!\!\int_{0}^{2\pi}\!\! Y_{lm}(\theta,\phi) \hspace*{1pt} Y_{l'\!,-m}(\theta,\phi) \hspace*{1pt} &Y_{l''\!,-m''}(\theta,\phi) \hspace*{2pt} \mathrm{d}\phi \hspace*{1pt} \sin\theta \hspace*{2pt} \mathrm{d}\theta\\\notag
\hspace*{2pt}&=\hspace*{2pt} \sqrt{\frac{(2l+1) (2l'\!+1) (2l''\!+1)}{4\pi}}
\begin{pmatrix}
l \!&\! l' \!&\! l'' \\
0 \!&\! 0  \!&\! 0
\end{pmatrix}\!\!
\begin{pmatrix}
l \!&\!   l'  \!&\!   l''  \\
m \!&\!\! -m \!&\!\! -m''
\end{pmatrix}\!.
\end{align}

The second Wigner-$3j$ symbol on the right-hand side of \eqref{eq:Gaunt} vanishes unless $l+l'\!+l''$ is even and $m''\!= 0$ \citep[Eqs.\ 3.177 and 3.195]{biedenharn_louck}. We can now insert
\begin{equation*}
Y_{l''\!,0}(0,0) \hspace*{2pt}=\hspace*{2pt} \sqrt{\frac{2l''\!+1}{4\pi}}
\end{equation*}
into \eqref{eq:spherical_bessel_addition_theorem_zwischenschritt_2}.
Relabeling $l''\! = k$ and using the triangle rule for the Wigner-$3j$ symbols \citep[Eq.\:3.191]{biedenharn_louck} gives the final result \eqref{eq:Bessel_addition_theorem}. The uniform convergence with respect to $[\vartheta'\!,\varphi]$ follows again from \eqref{eq:Bessel_estimate} and \eqref{eq:Y_estimate} via Weierstrass' test.\hfill\qedsymbol

\begin{remark}\label{rem:negative_z}
When considering translations in the \emph{negative} direction of the $z$ axis (\textit{i.e.}, switching to $-\nu < 0$), the additional factor $(-1)^{k}$ appears on the right-hand side of \eqref{eq:Bessel_addition_theorem}. This can easily be seen by substituting $\gamma + \pi$ for $\gamma$ in the above proof and noting that (cf.\ \citep[Sec.\:1.6.2]{dai_xu})
\begin{equation*}
Y_{lm}(\vartheta + \pi,\varphi) \hspace*{2pt}=\hspace*{2pt} (-1)^{l+m} \hspace*{1pt} Y_{l,m}(\vartheta,\varphi), ~~~~ |m| \leq l \in \mathbb{N}_{0}.
\end{equation*}
\end{remark}

For later use, we note the following property of the coefficients $A_{k}^{(ll'|m|)}$:

\begin{lemma}\label{lemma:A_uneven}
If \hspace*{1pt}$l-l'\!+k$ is odd, then $A_{k}^{(ll'|m|)} = 0$.
\end{lemma}

\noindent\textit{Proof.} This follows from the fact that $l-l'\!+k = (l+l'\!+k)-2l'$ is odd if and only if $l+l'\!+k$ is odd, and in this case \citep[Eq.\:3.195]{biedenharn_louck}
\begin{flalign*}
&& 
\begin{pmatrix}
l \!&\! l' \!&\! k \\
0 \!&\! 0  \!&\! 0
\end{pmatrix}\!
\hspace*{2pt}=\hspace*{2pt}0.
&& \mathllap{\qedsymbol}
\end{flalign*}

We now introduce the tools that we use to establish a connection between the functions $j_{l}(r) Y_{lm}(\vartheta,\varphi)$ and the SGL basis functions. In the GTO case, this was done via a certain class of Hankel transforms, which the author referred to as the spherical Bessel transform \citep[Sec.\:2.3]{ritchie}. The fact that the GTO basis functions are eigenfunctions of this transform simplifies the derivation of the GTO translation matrix elements crucially. As opposed to this, due to divergence of the underlying integrals, the spherical Bessel transform is not even well-defined in the SGL case. Therefore, we introduce a Gauss-Weierstrass convergence factor here:

\begin{definitioncited}[\textbf{(Weighted spherical Bessel transform)}]\label{def:weighted_spherical_Bessel_transform}
Let $l \in \mathbb{N}_{0}$ and $\gamma > 0$ be given. For admissible functions $f : (0,\infty) \to \mathbb{R}$, we define the \textit{weighted spherical Bessel transform} as
\begin{equation*}
\tilde{f}_{l}(\gamma; \beta) \hspace*{2pt}\coloneqq\hspace*{2pt} \sqrt{\frac{2}{\pi}} \int_{0}^{\infty}\!\! f(\xi) \hspace*{1pt} j_{l}(\beta \xi) \hspace*{1pt} \xi^{2} \hspace*{1pt} \mathrm{e}^{-\gamma \xi^{2}} \mathrm{d}\xi, ~~~~ \beta \in (0,\infty).
\end{equation*}
\end{definitioncited}

For a large class of functions, including the radial part of the SGL basis functions, a corresponding inversion formula can be derived (see \citep[Sec.\:14.4]{watson}). This inversion formula allows for a pointwise reconstruction of such functions from their spherical Bessel transform.

\begin{lemmacited}[\textbf{(Inversion formula)}]\label{lemma:inversion_formula}
Let \hspace*{1pt}$l \in \mathbb{N}_{0}$ and \hspace*{1pt}$\gamma > 0$ be given. Let further $f : (0,\infty) \to \mathbb{R}$ be continuous and of bounded variation on every bounded subset of \hspace*{1pt}$(0,\infty)$, and let $f$ satisfy
\begin{equation*}
\int_{0}^{\infty}\!\! |f(\xi)| \hspace*{1pt} \xi \hspace*{1pt} \mathrm{e}^{-\gamma \xi^{2}} \mathrm{d}\xi \hspace*{2pt}<\hspace*{2pt} \infty.
\end{equation*}
Then $f$ can be recovered pointwise from its Bessel transform via
\begin{equation*}\label{eq:Bessel_inverse}
f(\xi) \hspace*{2pt}=\hspace*{2pt} \mathrm{e}^{\gamma \xi^{2}} \sqrt{\frac{2}{\pi}} \int_{0}^{\infty}\!\! \tilde{f}_{l}(\gamma; \beta) \hspace*{1pt} j_{l}(\beta \xi) \hspace*{1pt} \beta^{2} \hspace*{1pt} \mathrm{d}\beta, ~~~~ \xi \in (0,\infty).
\end{equation*}
\end{lemmacited}

In the derivation of the SGL translation matrix elements in the upcoming Section \ref{sec:proof}, technical difficulties arise from the introduction of the Gauss-Weierstrass weight to the spherical Bessel transform. Specifically, when solving these problems, we encounter the hypergeometric series:

\begin{definitioncited}[\textbf{(Hypergeometric series)}]\label{def:hypergeometric_series}
For $p,q\in\mathbb{N}$ with $p \leq q + 1$ and parameters $a_{0},...,a_{p-1} \in \mathbb{C}$; $b_{0},...,b_{q-1} \in \mathbb{C} \setminus \lbrace0,-1,-2,\dots\rbrace$, the \emph{hypergeometric series} ${}_{p}F_{q}(a_{0},...,a_{p-1};$ $b_{0},...,b_{q-1}) : \mathbb{C} \to \mathbb{C}$ is defined as the power series \citep[Eq.\:2.1.2]{andrews_askey_roy}
\begin{equation}\label{eq:hypergeometric_series}
{}_{p}F_{q}(a_{0},\dots,a_{p-1}; b_{0},\dots,b_{q-1}; z) \hspace*{2pt}\coloneqq\hspace*{2pt} \sum_{k=0}^{\infty} \frac{(a_{0})_{k} \cdots (a_{p-1})_{k}}{(b_{0})_{k} \cdots (b_{q-1})_{k}} \frac{z^{p}}{p!}, ~~~~ z \in \mathbb{C},
\end{equation}
with the \emph{Pochhammer symbol (rising factorial)}
\begin{equation*}
(c)_{k} \hspace*{1pt}\coloneqq\hspace*{2pt} c(c+1)(c+2)\cdots(c+k-1), ~~~~ c \in \mathbb{C}, ~ k \in \mathbb{N}_{0}.
\end{equation*}
\end{definitioncited}

\begin{remark}
By definition, the order of the parameters $a_{0},...,a_{p-1}$ and $b_{0},...,b_{q-1}$ is irrelevant. If one of the parameters $a_{0},...,a_{p-1}$ is a non-positive integer, then by definition of the Pochhammer symbol, the hypergeometric series reduces to a polynomial.
\end{remark}

In this work, we will only encounter the special cases ${}_{1}F_{1}$ and ${}_{2}F_{1}$. The series ${}_{1}F_{1}$ is called \emph{Kummer's (confluent hypergeometric) function}. It is absolutely and uniformly convergent on bounded subsets of the complex plane \citep[Thm.\:2.1.1]{andrews_askey_roy}. The function ${}_{2}F_{1}$, commonly referred to as \emph{the} (\emph{ordinary} or \emph{Gaussian}) \emph{hypergeometric function}, will appear only as a polynomial in this work. 

Finally, we state two auxiliary lemmas for later use. The first lemma is well-known in the theory of finite difference equations. 

\begin{lemma}\label{lemma:finite_differenzen}
Let \hspace*{1pt}$p$ be a polynomial of degree less than $n \in \mathbb{N}$. Then
\begin{equation*}
\sum_{k=0}^{n} (-1)^{k} \binom{n}{k} \hspace*{1pt} p(k) \hspace*{2pt}=\hspace*{2pt} 0.
\end{equation*}
\end{lemma}

\begin{prf}
See \citep[Sec.\:1.4]{levy_lessman}, for instance.
\end{prf}

\begin{lemma}\label{lemma:gamma_trick}
Let \hspace*{1pt}$n\in\mathbb{N}_{0}$. Then
\begin{equation*}
\Gamma(n+1/2) \hspace*{2pt}=\hspace*{2pt} \sqrt{\pi} \hspace*{1pt} (n+1)_n \hspace*{1pt} 4^{-n}.
\end{equation*}
\end{lemma}

\begin{prf}
This follows from the functional equation $\Gamma(z+1) = z\Gamma(z)$, $z\in\mathbb{C}\setminus\lbrace0,-1,-2,\dots\rbrace$ 
\citep[Eq.\:1.1.6]{andrews_askey_roy} and the fact that $\Gamma(1/2) = \sqrt{\pi}$ 
\citep[Eq.\:1.1.22]{andrews_askey_roy}. 
\end{prf}

\section{SGL translation matrix elements}\label{sec:proof}
We now derive the closed-form expression for the SGL translation matrix elements. To this end, let $|m| \leq l < n \in \mathbb{N}$ and $|m'| \leq l'\! < n'\! \in \mathbb{N}$ be given. Moreover, fix $\nu > 0$ and set $\boldsymbol{x}'\! \coloneqq \boldsymbol{x} - \nu \mathbf{e}_{z} = [r'\!,\vartheta'\!,\varphi]$ for $\boldsymbol{x} = [r,\vartheta,\varphi]$. With the dominated convergence theorem, 
we find that
\begin{equation*}
T_{nn'll'}^{(mm')}(\nu)
\hspace*{2pt}=\hspace*{2pt} \lim_{\gamma \to 0+} \hspace*{1pt} \int_{0}^{\infty} \!\!\!\int_{0}^{\pi} \!\!\int_{0}^{2\pi}\!\! R_{nl}(r') \hspace*{0pt} Y_{lm}(\vartheta'\!,\varphi) \hspace*{1pt} \mathrm{e}^{-\gamma (r')^{2}} R_{n'l'}(r) \hspace*{0pt} \overline{Y_{l'm'}(\vartheta,\varphi)} \hspace*{1pt} \mathrm{e}^{-r^{2}} \mathrm{d}\varphi \hspace*{1pt} \sin\vartheta \hspace*{2pt} \mathrm{d}\vartheta \hspace*{2pt} r^{2} \hspace*{1pt} \mathrm{d}r.
\end{equation*}
By Lemma \ref{lemma:inversion_formula}, the first part of Lemma \ref{lemma:Bessel_addition_theorem}, and with Remark \ref{rem:negative_z}, we get for $r' > 0$
\begin{align}\label{eq:steppaaa}
R_{nl}(r') Y_{lm}(\vartheta'\!,\varphi) \mathrm{e}^{-\gamma (r')^{2}\!}\!
\hspace*{0pt}&=\hspace*{0pt} \!\sqrt{\frac{2}{\pi}}\!\hspace*{0pt} \int_{0}^{\infty}\!\!\! \hspace*{0.5pt}\tilde{R}_{nl}(\gamma;\beta) \hspace*{0pt} j_{l}(\beta r') \hspace*{0pt} Y_{lm}(\vartheta'\!,\varphi) \hspace*{0pt} \beta^{2} \hspace*{0pt} \mathrm{d}\beta\\[2pt]
\hspace*{0pt}&=\hspace*{0pt} \!\sqrt{\frac{2}{\pi}}\!\hspace*{0pt} \int_{0}^{\infty}\!\!\! \hspace*{0.5pt}\tilde{R}_{nl}(\gamma;\beta) \!\sum_{l''\!=0}^{\infty} \hspace*{0pt} \sum_{k = |l - l''|}^{l + l''\!}\!\!\! (-1)^{k} A_{k}^{(ll''|m|)\!} j_{k}(\beta \nu) j_{l''}(\beta r) Y_{l''m}(\vartheta,\varphi) \beta^{2} \hspace*{0pt} \mathrm{d}\beta.\notag
\end{align}
It is \citep[Sec.\:8.9, Eq.\:5]{erdelyi}
\begin{equation}\label{eq:bessel_transform_gamma}
\tilde{R}_{nl}(\gamma;\beta) \hspace*{2pt}=\hspace*{2pt} N_{nl} \hspace*{1pt} \frac{(\gamma-1)^{n-l-1}}{\gamma^{n+1/2}} \hspace*{1pt} \beta^{l} L_{n-l-1}^{(l+1/2)}\bigg(\frac{\beta^{2}}{4 \gamma (1-\gamma)}\bigg) \hspace*{1pt} \mathrm{e}^{-\beta^{2}/4\gamma} \hspace*{1pt} \bigg(\frac{1}{2}\bigg)^{\!l+3/2}, ~~~~ 0 < \gamma < 1.
\end{equation}
%
%
For $\gamma = 1$, the formula reads as \citep[Sec.\:8.9, Eq.\:2]{erdelyi}
\begin{equation}\label{eq:bessel_transform_1}
\tilde{R}_{nl}(1;\beta) \hspace*{2pt}=\hspace*{2pt} \frac{N_{nl}}{(n-l-1)!} \hspace*{1pt} \beta^{2n-l-2} \hspace*{1pt} \mathrm{e}^{-\beta^{2}/4} \hspace*{1pt} \bigg(\frac{1}{2}\bigg)^{\!2n-l-1/2}.
\end{equation}
Formula \eqref{eq:bessel_transform_gamma} shows us particularly that after inserting \eqref{eq:steppaaa} into the above formula for $T_{nn'll'}^{(mm')}(\nu)$, changing the order of integration is legitimate. By Fubini's theorem 
and the second part of Lemma \ref{lemma:Bessel_addition_theorem}, we thus obtain 
\begin{align}\label{eq:zwischenschritt}
T_{nn'll'}^{(mm')}(\nu)
\hspace*{0pt}&=\hspace*{0pt} \lim_{\gamma \to 0+} \sqrt{\frac{2}{\pi}} \int_{0}^{\infty}\!\! \int_{0}^{\infty}\!\! \tilde{R}_{nl}(\gamma;\beta) \sum_{l''=0}^{\infty} \sum_{k = |l - l''|}^{l + l''\!}\!\!\! (-1)^{k} A_{k}^{(ll''|m|)} \delta_{l'l''} \hspace*{1.5pt} \delta_{mm'} \\ 
&\hphantom{= \lim_{\gamma \to 0+}}\hspace*{14pt} \times \hspace*{1pt} j_{k}(\beta\nu) \hspace*{1pt} j_{l''}(\beta r) \hspace*{1pt} \beta^{2} \hspace*{1pt} \mathrm{d}\beta \hspace*{1pt} R_{n'l'}(r) \hspace*{1pt} r^{2} \hspace*{1pt} \mathrm{e}^{-r^{2}} \mathrm{d}r \notag\\[2pt]
%
%
\hspace*{0pt}&=\hspace*{0pt} \delta_{mm'\!}\!\! \sum_{k = |l - l'|}^{l + l'\!}\!\!\! (-1)^{k} A_{k}^{(ll'|m|)\!} \lim_{\gamma \to 0+} \underbrace{\int_{0}^{\infty}\!\! \tilde{R}_{nl}(\gamma;\beta) \hspace*{1pt} \tilde{R}_{n'l'}(1;\beta) \hspace*{1pt} j_{k}(\beta\nu) \hspace*{1pt} \beta^{2} \hspace*{1pt} \mathrm{d}\beta}_{\hphantom{I^{(nn'll')}_{k}(\gamma;\nu)} \hspace*{3.5pt} \eqqcolon \hspace*{2pt} I^{(nn'll')}_{k}(\gamma;\nu)},\notag
\end{align}
where the $\lim_{\gamma \to 0+} I^{(nn'll')}_{k}(\gamma;\nu)$ is yet to be investigated. The two Kronecker symbols $\delta_{l'l''}$ and $\delta_{mm'}$ are due to the orthonormality of the spherical harmonics. Because of the remaining Kronecker symbol $\delta_{mm'}$, we shall write $T_{nn'll'}^{(mm')}(\nu) = T_{nn'll'}^{(|m|)}(\nu)$. Due to Lemma \ref{lemma:A_uneven}, we have to consider only the case when $l-l'\!+k$ is even in the following.

\begin{lemma}\label{lemma:lim_I}
Let \hspace*{1pt}$\hspace*{1pt}l-l'\!+k$ be even and set \hspace*{1pt}$\mu \coloneqq n'\!+(l-l'\!+k)/2 \in \mathbb{N}$. Then
\begin{align*}
\lim_{\gamma \to 0+} I^{(nn'll')}_{k}(\gamma;\nu) \hspace*{2pt}&=\hspace*{2pt} (-1)^{n-l-1} \hspace*{1pt} \nu^{k} \hspace*{1pt} \frac{\sqrt{\pi}}{4} \hspace*{1pt} \frac{N_{nl} \hspace*{1pt} N_{n'l'}}{(n'\!-l'\!-1)!} \\
&\hspace*{16pt}\times \sum_{j=0}^{n-l-1} \frac{(-1)^{j}}{j!} \binom{n-1/2}{n-l-1-j} \frac{\Gamma(\mu+j+1/2)}{\Gamma(k+3/2)} \hspace*{1pt} C_{j}^{(nn'll'k)}(\nu),
\end{align*}
where
\begin{equation*}
C_{j}^{(nn'll'k)}(\nu) \hspace*{2pt}=\hspace*{2pt}
\sum\limits_{p=0}^{n-\mu} \frac{{}_{2}F_{1}(j,-n+\mu+p;-n+\mu-j+1;-1) (-\mu+k-j+1)_{p}}{p! \hspace*{1pt} (k+3/2)_{p}} \binom{-j-p}{n-\mu-p} \hspace*{1pt} \nu^{2p}.
\end{equation*}
\end{lemma}

\noindent\textit{Proof.} By inserting \eqref{eq:bessel_transform_gamma} and \eqref{eq:bessel_transform_1} into $I^{(nn'll')}_{k}(\gamma;\nu)$, we derive
\begin{align*}\label{eq:zwischenschritt_2}
I^{(nn'll')}_{k}(\gamma;\nu)
\hspace*{2pt}&=\hspace*{2pt} (-1)^{n-l-1} \frac{N_{nl} \hspace*{1pt} N_{n'l'}}{(n'\!-l'\!-1)!}\! \sum_{j=0}^{n-l-1} \frac{(-1)^{j}}{j!} \binom{n-1/2}{n-l-1-j} \! \bigg(\frac{1}{2}\bigg)^{\!2n'\!+l-l'\!+2j+1}\\
&\hspace*{16pt} \times \gamma^{-n-j-1/2} \hspace*{1.5pt} (1-\gamma)^{n-l-j-1} \underbrace{\int_{0}^{\infty}\!\! \mathrm{e}^{-\beta^{2}(1+\gamma^{-1})/4} \hspace*{1pt} j_{k}(\beta\nu) \hspace*{1pt} \beta^{2n'\!+l-l'+2j} \hspace*{1pt} \mathrm{d}\beta
}_{\hphantom{J^{(n'll'k)}_{j}(\gamma;\nu)} \hspace*{3.5pt} \eqqcolon \hspace*{2pt} J^{(n'll'k)}_{j}(\gamma;\nu)},
\end{align*}
where we have also used the closed-form expression of the generalized Laguerre polynomials \citep[Eq.\:6.2.2]{andrews_askey_roy}.
Now we solve the integral on the right-hand side by using [Erd\'{e}lyi, 1954, Sec.\:8.6, Eq.\:14]:
\begin{align*}
J^{(n'll'k)}_{j}(\gamma;\nu) \hspace*{2pt}=\hspace*{4.5pt} &\sqrt{\frac{\pi}{2}} \hspace*{2pt} \frac{\Gamma(\mu+j+1/2)}{\Gamma(k+3/2)} \bigg(\frac{\gamma}{1+\gamma}\bigg)^{\!\mu+j+1/2\!} \hspace*{1pt} 2^{2n'\!+l-l'\!+2j-1/2\!} \\ 
&\hspace*{9.5pt}\times \nu^{k} \hspace*{1.5pt} \mathrm{e}^{- \frac{\gamma}{1+\gamma} \nu^{2}} \hspace*{1.5pt} {}_{1}F_{1} \bigg(\!\!-\mu+k-j+1; \hspace*{1pt} k+3/2; \hspace*{1pt} \frac{\gamma}{1+\gamma} \nu^{2}\bigg),
\end{align*}
where ${}_{1}F_{1}$ denotes Kummer's confluent hypergeometric function (cf.\ Def.\:\ref{def:hypergeometric_series}).

Now we investigate the limit of
\begin{align*}
I_{k}^{(nn'll')}(\gamma;\nu)
\hspace*{2pt}&=\hspace*{2pt} (-1)^{n-l-1} \hspace*{1pt} \nu^{k} \hspace*{1pt} \frac{\sqrt{\pi}}{4} \hspace*{1pt} \frac{N_{nl} \hspace*{1pt} N_{n'l'}}{(n'\!-l'\!-1)!} \hspace*{1pt} \gamma^{\mu-n} \hspace*{1pt} \mathrm{e}^{- \frac{\gamma}{1+\gamma} \nu^{2}} \\
&\hspace*{16pt}\times \sum_{j=0}^{n-l-1} \frac{(-1)^{j}}{j!} \binom{n-1/2}{n-l-1-j} \frac{\Gamma(\mu+j+1/2)}{\Gamma(k+3/2)} \hspace*{1pt} (1+\gamma)^{-\mu-j-1/2}\\
&\hspace*{16pt}\times (1-\gamma)^{n-l-j-1} {}_{1}F_{1}\bigg(\!\!-\mu+k-j+1; \hspace*{1pt} k+3/2; \hspace*{1pt} \frac{\gamma}{1+\gamma} \nu^{2}\bigg)
\end{align*}
as $\gamma$ tends to zero from above.
Since the right-hand side contains the critical factor $\gamma^{\mu-n}$, we distinguish between three different cases.

If $k > 2(n-n')-l+l'\!$, then $\mu-n > 0$, and thus $\lim_{\gamma \to 0+} I_{k}^{(nn'll')}(\gamma;\nu) = 0$, since by \eqref{eq:hypergeometric_series},
\begin{equation}\label{eq:M_limit}
\lim_{\gamma \to 0+} {}_{1}F_{1}\bigg(\!\!-\mu+k-j+1; \hspace*{1pt} k+3/2; \hspace*{1pt} \frac{\gamma}{1+\gamma} \nu^{2}\bigg) \hspace*{2pt}=\hspace*{2pt} 1.
\end{equation}
If, on the other hand, $k = 2(n-n')-l+l'\!$, then $\mu-n = 0$, and
\begin{align*}
\lim_{\gamma \to 0+} I_{k}^{(nn'll')}(\gamma;\nu) 
\hspace*{2pt}&=\hspace*{2pt} (-1)^{n-l-1} \hspace*{1pt} \nu^{k} \hspace*{1pt} \frac{\sqrt{\pi}}{4} \hspace*{1pt} \frac{N_{nl} \hspace*{1pt} N_{n'l'}}{(n'\!-l'\!-1)!} \\
&\hspace*{16pt}\times \sum_{j=0}^{n-l-1} \frac{(-1)^{j}}{j!} \binom{n-1/2}{n-l-1-j} \frac{\Gamma(\mu+j+1/2)}{\Gamma(k+3/2)},
\end{align*}
which also follows from \eqref{eq:M_limit}.

If $k < 2(n-n')-l+l'\!$, then $\mu-n < 0$. In this case, we cannot make use of \eqref{eq:M_limit} directly. 
However, we have the power-series representation
\begin{equation*}
(1+\gamma)^{-j-p} \hspace*{1pt} (1-\gamma)^{-j}  \hspace*{2pt}=\hspace*{2pt} \sum_{q=0}^{\infty} \binom{-j-p}{q} \hspace*{1pt} {}_{2}F_{1}(j,-q;1-q-p-j;-1) \hspace*{1pt} \gamma^{q},
\end{equation*}
with the hypergeometric function ${}_{2}F_{1}$ (cf.\ Def.\:\ref{def:hypergeometric_series}), derived from the Cauchy product 
of the binomial series of $(1+\gamma)^{-j-p}$ and $(1-\gamma)^{-j\!}$ \citep[Eq.\:3.6.8]{abramowitz_stegun}. This series is absolutely convergent for $|\gamma| < 1$.
Rearranging the summands, we thus obtain together with \eqref{eq:hypergeometric_series} for sufficiently small $0 < \gamma < 1$
\begin{align}\label{eq:zwischenschritt_3}
I_{k}^{(nn'll')}(\gamma;\nu) 
\hspace*{2pt}&=\hspace*{2pt} (-1)^{n-l-1} \hspace*{1pt} \nu^{k} \hspace*{1pt} \frac{\sqrt{\pi}}{4} \hspace*{1pt} \frac{N_{nl} \hspace*{1pt} N_{n'l'}}{(n'\!-l'\!-1)!} \hspace*{1pt} \mathrm{e}^{- \frac{\gamma}{1+\gamma} \nu^{2}} \hspace*{1pt} \Gamma(k+3/2)^{-1}\\
&\hspace*{16pt}\times (1+\gamma)^{-\mu-1/2} \hspace*{1pt} (1-\gamma)^{n-l-1} \sum_{p,q=0}^{\infty} \frac{\nu^{2p}}{p! \hspace*{1pt} (k+3/2)_{p}} \hspace*{1pt} D_{pq}^{(nn'll'k)} \hspace*{1pt} \gamma^{p+q-n+\mu},\notag
\end{align}
where we set
\begin{align}\label{eq:def_D}
D_{pq}^{(nn'll'k)} \hspace*{2pt}&\coloneqq\hspace*{2pt} \sum_{j=0}^{n-l-1} \frac{(-1)^{j}}{j!} \binom{n-1/2}{n-l-1-j} \Gamma(\mu+j+1/2)  \\
&\hspace*{16pt}\times (-\mu+k-j+1)_{p} \hspace*{1pt} \binom{-j-p}{q} \hspace*{1pt} {}_{2}F_{1}(j,-q;1-q-p-j;-1).\notag
\end{align}

Before taking the limit $\gamma \to 0+$, we now prove that $D_{pq}^{(nn'll'k)} = 0$ for $p+q < n-\mu$. To this end, we first use Lemma \ref{lemma:gamma_trick} to rewrite
\begin{align*}
\binom{n-1/2}{n-l-1-j}\! 
\hspace*{2pt}&=\hspace*{2pt} \frac{\Gamma(n+1/2)}{(n-l-1-j)! \hspace*{1pt} \Gamma(l+j+3/2)}\label{eq:binom_1}\\[2pt]
\hspace*{2pt}&=\hspace*{2pt}  \hspace*{0pt} \frac{4^{-n+l+j+1} \hspace*{1pt} (n+1)_{n}}{(n-l-1-j)! \hspace*{1pt} (l+j+2)_{l+j+1}}.\notag
\end{align*}
Again by Lemma \ref{lemma:gamma_trick}, we find that
\begin{equation*}
\Gamma(\mu+j+1/2) \hspace*{2pt}=\hspace*{2pt} \sqrt{\pi} \hspace*{1pt} (\mu+j+1)_{\mu+j} \hspace*{1pt} 4^{-\mu-j}.
\end{equation*}
The fact that 
\begin{equation*}
\mu \hspace*{2pt}\geq\hspace*{2pt} 1 + l'\! + \frac{l - l'\! + k}{2} \hspace*{2pt}\geq\hspace*{2pt} 1 + \frac{l + l'\! + |l - l'|}{2} \hspace*{2pt}\geq\hspace*{2pt} 1 + l
\end{equation*}
then reveals
\begin{align}
\binom{n-1/2}{n-l-1-j} \Gamma(\mu+j+1/2)\label{eq:komponente_1}
\hspace*{2pt}&=\hspace*{2pt} \frac{\sqrt{\pi} \hspace*{1pt} (n+1)_{n}}{(n-l-1-j)!} \frac{(\mu+j+1)_{\mu+j}}{(l+j+2)_{l+j+1}} \hspace*{1pt} 4^{\mu-n+l+1} \\[2pt]
\hspace*{2pt}&=\hspace*{2pt} \frac{\sqrt{\pi} \hspace*{1pt} (n+1)_{n}}{(n-l-1-j)!} \hspace*{1pt} (l+j+3/2)_{\mu-l-1} \hspace*{1pt} 4^{2\mu-n}.\notag
\end{align}
By definition of the hypergeometric function ${}_{2}F_{1}$ (Def.\:\ref{def:hypergeometric_series}), we get
\begin{align}
\binom{-j-p}{q} \hspace*{1pt} {}_{2}F_{1}(j,-q;1-j-q-p;-1)\label{eq:komponente_2}
\hspace*{2pt}&=\hspace*{2pt} \frac{1}{q!} \sum_{s=0}^{q} \hspace*{1pt} (-1)^{s} \hspace*{1pt} (-q)_{s} \hspace*{1pt} (1-p-q+s-j)_{q-s} \hspace*{1pt} (j)_{s}.
\end{align}
Inserting \eqref{eq:komponente_1} and \eqref{eq:komponente_2} into \eqref{eq:def_D}, changing the order of summation, and expanding with $(n-l-1)!$ yields
\begin{align*}
D_{pq}^{(nn'll'k)\!} \hspace*{2pt}&=\hspace*{2pt} 4^{2\mu-n} \hspace*{1pt} \frac{\sqrt{\pi} \hspace*{1pt} (n+1)_{n}}{(n-l-1)!} \hspace*{1pt} \frac{1}{q!} \sum_{s=0}^{q} \hspace*{1pt} (-1)^{s} \hspace*{1pt} (-q)_s \\
&\hspace*{16pt}\times \sum_{j=0}^{n-l-1}\! (-1)^{j} \binom{n\!-\!l\!-\!1}{j} (l\!+\!j\!+\!3/2)_{\mu-l-1} \hspace*{1pt} (\!-\mu\!+\!k\!-\!j\!+\!1)_{p} \hspace*{1pt} (1\!-\!p\!-\!q\!+\!s\!-\!j)_{q-s} \hspace*{1pt} (j)_{s}.
\end{align*}
Since $(l+j+3/2)_{\mu-l-1} \hspace*{1pt} (-\mu+k-j+1)_{p} \hspace*{1pt} (1-p-q+s-j)_{q-s} \hspace*{1pt} (j)_{s}$ can be seen as a polynomial of degree $p+q+\mu-l-1$ evaluated at $j$, we find that $D_{pq}^{(nn'll'k)\!} = 0$ for $p + q < n - \mu$ by Lemma \ref{lemma:finite_differenzen}, as claimed. Therefore, the double series in \eqref{eq:zwischenschritt_3} can be seen as a power series in $\gamma$ with convergence radius one, which as such converges uniformly on all compact subsets of $[0,1)$. This allows for a term-by-term limit formation $\gamma \to 0+$.
Because all summands in \eqref{eq:zwischenschritt_3} with $p+q > n-\mu$ vanish, we are left with the summands where
$p+q = n-\mu$. This means that
\begin{flalign*}
&&
\lim_{\gamma \to 0+} I_{k}^{(nn'll')}(\gamma;\nu) 
\hspace*{2pt}=\hspace*{2pt} (-1)^{n-l-1} \hspace*{1pt} \nu^{k} \hspace*{1pt} \frac{\sqrt{\pi}}{4} \hspace*{1pt} \frac{N_{nl} \hspace*{1pt} N_{n'l'}}{(n'\!-l'\!-1)!} \hspace*{1pt}\sum\limits_{p=0}^{n-\mu}\hspace*{1pt} \frac{\nu^{2p}}{p! \hspace*{1pt} (k+3/2)_{p}} \hspace*{1pt} D_{p,(n-\mu)-p}^{(nn'll'k)}.
&& \mathllap{\qedsymbol}
\end{flalign*}

In summary, applying Lemma \ref{lemma:lim_I} to \eqref{eq:zwischenschritt}, we obtain the following closed-form expression for the SGL translation matrix elements:

\begin{maintheoremcited}[\textbf{(SGL translation matrix elements)}]\label{thm:sgl_translation_matrix_elements}
The SGL translation matrix elements possess the closed-form expression
\begin{align*}
T_{nn'll'}^{(|m|)}(\nu)
\hspace*{2pt}&=\hspace*{2pt} (-1)^{n-l-1} \hspace*{1pt} \frac{\sqrt{\pi}}{4} \hspace*{1pt} \frac{N_{nl} \hspace*{1pt} N_{n'l'}}{(n'\!-l'\!-1)!} \sum_{k = |l-l'|}^{l+l'\!} (-1)^{k} A_{k}^{(ll'|m|)} \hspace*{1pt} \nu^{k} \\
&\hspace*{16pt}\times \sum_{j=0}^{n-l-1} \frac{(-1)^{j}}{j!} \binom{n-1/2}{n-l-1-j} \frac{\Gamma(\mu+j+1/2)}{\Gamma(k+3/2)} \hspace*{1pt} C_{j}^{(nn'll'k)}(\nu), ~~~~ \nu > 0,
\end{align*}
where the coefficients $A_{k}^{(ll'|m|)}$ are those of Lemma \ref{lemma:Bessel_addition_theorem}, while $\mu$ and the functions $C_{j}^{(nn'll'k)}$ are given in Lemma \ref{lemma:lim_I} (setting \hspace*{1pt}$C_{j}^{(nn'll'k)} \coloneqq 0$ whenever $l-l'\!+k$ is odd).
\end{maintheoremcited}

\begin{remark}
Note that the hypergeometric function ${}_{2}F_{1}$ in $C_{j}^{(nn'll'k)}$ is, in fact, a polynomial.
\end{remark}

\section*{Acknowledgments}
The authors would like to thank Sabrina Kombrink, Vitalii Myroniuk, and Nadiia Derevianko for scientific discussion and helpful comments on the manuscript.

\setlength{\bibsep}{3pt}
\bibliography{references}
\bibliographystyle{my_apa}

\end{document}